\input amstex
\documentstyle{amsppt}
\topmatter
\title
Regularity of the Moduli Space of Instanton
Bundles $MI_{\Bbb P^3}(5)$.
\endtitle
\author Pavel I. Katsylo$^*$ and Giorgio Ottaviani$^{**}$
\endauthor
\thanks
$^*$ Supported by CRDF, grant RM1--206, and
INTAS, grant INTAS--OPEN--97--1570.
$^{**}$ partially supported by MURST funds.
\endthanks
\address
Moscow Independent University, Bol'shoi Vlas'evski\v i per., 11,
117463, Mo\-scow, Russia
\endaddress
\email
katsylo\@katsylo.mccme.ru
\endemail
\address
Dipartimento di Matematica, viale Morgagni 67 A, 50134 Firenze, Italy
\endaddress
\email
ottavian\@math.unifi.it
\endemail
\endtopmatter
\document

\def\kappa{\varkappa}

\def\hsigma{\widehat\sigma}
\def\tsigma{\widetilde\sigma}
\def\rig#1{\smash{
   \mathop{\longrightarrow}\limits\sp {#1}}}
\def\dow#1{\Big\downarrow
   \rlap{$\vcenter{\hbox{$\scriptstyle#1$}}$}}

\define\C{\Bbb C}
\define\PP{\Bbb P}

\define\SL{\operatorname{SL}}
\define\SP{\operatorname{Sp}}

\define\EE{\operatorname{E}}
\define\Ker{\operatorname{Ker}}
\define\rk{\operatorname{rk}}
\define\row{\operatorname{row}}
\def\Im{\operatorname{Im}}
\def\dim{\operatorname{dim}}
\define\Coker{\operatorname{Coker}}

\define\OO{\Cal O}

\define\qqq{\hfill $\square$}


\head \S0. Introduction. \endhead

Instanton vector bundles were defined by Atiyah, Drinfeld, Hitchin and
Manin [ADHM] in order to construct all the self-dual solutions of the
Yang-Mills equation over $S^4$. A mathematical instanton bundle $E$ on $\PP^3:=
\PP^3(\C)$ can be defined as the cohomology bundle of a monad
$$
\OO(-1)^k \ \rig{} \ \OO^{2k+2} \ \rig{} \ \OO(1)^k
$$
on $\PP^3$, where $c_2(E)=k$. This is equivalent to the condition that $E$ is a 
stable bundle of rank $2$ on $\PP^3$ such that $c_1(E)=0$, $c_2(E)=k$ and
$H^1(E(-2))=0$. If $E$ is a mathematical instanton bundle, then it is easy to
check by using Hirzebruch-Riemann-Roch Theorem that $h^1(S^2E)-h^2(S^2E)=8 k-3$ 
and so $8k -3$ is the expected dimension
of the moduli space of mathematical instanton bundles $MI_{\PP^3}(k)$.
It is not known if the moduli space $MI_{\PP^3}(k)$
is a regular variety of pure dimension $8 k - 3$.
The answer is affirmative in the cases $1\le k\le 4$ and this was proved
in [H], [ES] and [LeP]. In this article we extend these results to the case
$k=5$. More precisely, we prove the following

\proclaim{Theorem 0.1} For $2 \le k \le 5$ the moduli space
$MI_{\PP^3}(k)$ of mathematical instantons is a regular variety
of pure dimension $8 k - 3$.
\endproclaim

Our result should be compared with [AO2], where it is proved that the closure
of  $MI_{\PP^3}(5)$ in the Maruyama scheme of vector bundles of rank $2$ with
$c_1=0$, $c_2=5$ contains singular points.
Our proof requires tools both from invariant theory and algebraic geometry.

In \S 1 we prove some algebraic lemmas. In \S 2 we give an invariant
theoretical description of $MI_{\PP^3}(k)$
and we prove some result about unstable planes.
Finally, in \S 3 we prove Theorem 0.1.

The authors are grateful to A. N. Tyurin for fruitful advices. The first author 
thanks the Italian GNSAGA-INDAM for financial support. The second author is 
a member of the European network EAGER.


\head \S 1. An invariant theoretical description of $MI_{\PP^3}(k)$.
\endhead

Our first goal is to describe the moduli space $MI_{\PP^3}(k)$  in terms of 
invariant theory and prove that for any $\in MI_{\PP^3}(k)$ and any plane $H$ in
$\PP^3$ $h^0(E_H)\le 1$.
The group $\SL_{2k+2}$ acts canonically on the space
$\C^{2k+2}$. Let $h_1, ... , h_{2k+2}$
be the standard basis of the space $\C^{2k+2}$ and
$h_1^\ast, ... , h_{2k+2}^\ast$ be the dual basis of the dual space
$\C^{(2k+2) \ast}$. Consider the 2-form
$$
\omega = \sum_{1 \le i \le k+1} h_i^\ast \wedge h_{k+1+i}^\ast
\ = \ \sum_{1 \le i \le k+1}
(h_i^\ast \otimes h_{k+1+i}^\ast - h_{k+1+i}^\ast \otimes h_i^\ast)
\ \in \ \wedge^2 \C^{(2k+2) \ast}.
$$
Let $\SP_{2k+2}$ be the stabilizer of the
2-form $\omega$ in the group $\SL_{2k+2}$. The 2-form $\omega$
defines canonically the $\SP_{2k+2}$-isomorphism
$\C^{(2k+2) \ast} \simeq \C^{2k+2}$. We have the canonical actions
of the group $\SL_4 \times \SL_k \times \SP_{2k+2}$ on the
spaces $\C^4, \ \C^{4 \ast}, \ \C^{2k+2}, \ \C^k$, $ \C^{k \ast},
\ \C^4 \otimes \C^k, ... $.

We have the quadratic $\left(\SL_4 \times \SL_k \times \SP_{2k+2}\right)$-morphism
$$
\gathered
\gamma \ : \ \C^{4 \ast} \otimes \C^{k \ast} \otimes \C^{2k+2}
\rightarrow
S^2 \C^{4 \ast} \otimes \wedge^2 \C^{k \ast}, \\
(f_i^\ast \otimes b_j^\ast \otimes h_l
\ f_{i'}^\ast \otimes b_{j'}^\ast \otimes h_{l'}) \ \mapsto
\ \frac{1}{2} \omega (h_l, h_{l'})
\ (f_i^\ast f_{i'}^\ast) \otimes (b_j^\ast \wedge b_{j'}^\ast)
\endgathered
$$
and the bilinear $\left(\SL_4 \times \SL_k \times \SP_{2k+2}\right)$-morphisms
$$
\gathered
\beta : \C^{4 \ast} \otimes \C^{k \ast} \otimes \C^{2k+2} \ \times
\ \C^{2k+2} \rightarrow \C^{4 \ast} \otimes \C^{k \ast}, \\
\beta (f_i^\ast \otimes b_j^\ast \otimes h_l,
\ h_{l'}) \ = \ \omega (h_l, h_{l'}) f_i^\ast \otimes b_j^\ast \\ 
\endgathered$$
and
$$  \gathered
\varepsilon \ : \ \C^{4 \ast} \otimes \C^{k \ast} \otimes \C^{2k+2}
\ \times \ \C^4 \otimes \C^k
\ \rightarrow \ \C^{2k+2}, \\
(f_i^\ast \otimes b_j^\ast \otimes h_l, \ f_{i'} \otimes b_{j'})
\ \mapsto \ \delta_{ii'} \delta_{jj'} h_l 
\endgathered
$$

Consider the following conditions for an element
$A \in \C^{4 \ast} \otimes \C^{k \ast} \otimes \C^{2k+2}$
\roster
\item"{$(E_1)$}" \ $\varepsilon (A, f \otimes b) \neq 0$ for all
$0 \neq f \in \C^4$, $0 \neq b \in \C^k$,
\item"$(E_2)$" \ $\gamma (A) = 0$,
\item"$(E_3)$" \ $\beta(A,h) \neq 0$ for all $0 \neq h \in \C^{2k+2}$.
\endroster

An element $A \in \C^{4 \ast} \otimes \C^{k \ast} \otimes \C^{2k+2}$
defines the sheaf morphism $\OO^{2k+2} \ \rig{A} \ \OO(1)^k$.
This morphism and the symplectic structure over
$\OO^{2k+2}$ define the sequence
$$
\OO(-1)^k \ \rig{f_A^\top } \ \OO^{2k+2} \ \rig{f_A} \ \OO(1)^k .\tag 1.1
$$
The condition $(E_1)$ means that $f_A$ is surjective or
that $\Ker~f_A$ is locally free.
The condition $(E_2)$ means that the above sequence is a complex. Therefore,
$(E_1)$ and $(E_2)$ together mean that it is a monad according to [BH].
The condition $(E_3)$ means moreover that the cohomology bundle $E$ of the monad
is a stable vector bundle.
It is well known (see e.g. [AO1], Theor. 2.8) that the conditions $(E_1)$
and $(E_2)$ imply $(E_3)$.

Set
$$
I_i = \{ A \in \C^{4 \ast} \otimes \C^{k \ast} \otimes \C^{2k+2}
\ | \ \text{the condition} \ (E_i) \ \text{holds for} \ A \}, $$
$$ I \ \colon = \ I_1 \cap I_2 \cap I_3
\ = \ I_1 \cap I_2.
$$

Recall that the quotient $I/G$,
where $G = \SL_k \times \SP_{2k+2} \times \C^\ast$ is isomorphic to
$MI_{\PP^3}(k)$ (see [BH]).
Moreover, the stabilizer of any $A \in I$
in the group $G$ is equal to the kernel $H$ of the action
of $G$ on $I$ and $\dim (H) = 0$. Therefore, any $G$-orbit
\ $G \cdot A$, \ where $A \in I$, is isomorphic to $G/H$
and $\dim (G \cdot A) \ = \ \dim (G) \ = \ 3k^2 + 5k + 3$.
Consider the canonical morphism
$$
\pi \ : \ I \ \rightarrow \ I/G \ = \ MI_{\PP^3}(k).
$$

The results above imply the following fact:
\proclaim{Lemma 1.1} For any $A \in I$ we have
$$
\dim (T_{\pi (A)}MI_{\PP^3}(k)) \ = \ \dim (T_A I) - 3k^2 - 5k - 3.
$$ 
\endproclaim

\proclaim {Theorem 1.2} Let $E$ be an instanton bundle on ${\PP}^3$ and
let $H$ be a plane. Then $h^0(E_{|H})\le 1$.
\endproclaim
\demo{Proof} The bundle $E$ is the cohomology bundle of a monad
$$
\OO(-1)^k \ \rig{} \ \OO^{2k+2} \ \rig{b} \ \OO(1)^k ,
$$
where $c_2(E)=k$.
Notice that it is enough to show that
$h^0(K_{|H})\le 1$, where $K:=\Ker b$. A section $s\in H^0(K_{|H})$ induces 
the diagram
$$
\matrix &&0&&0\\
&&\dow{}&&\dow{}\\
&&\OO_H&\rig{\simeq}&\OO_H\\
&&\dow{s}&&\dow{}\\
0&\rig{}&K_{|H}&\rig{i}&\OO_H^{2k+2}&\rig{}&\OO_H(1)^k&\rig{}&0\\
&&\dow{}&&\dow{}&&\dow{\simeq}\\
0&\rig{}&T&\rig{}&\OO_H^{2k+1}&\rig{}&\OO_H(1)^k&\rig{}&0\\
&&\dow{}&&\dow{}\\
&&0&&0 \endmatrix
$$
Since $i\circ s$ is injective, we see that $s$ splits $K_{|H}$
as $K_{|H}\simeq \OO_H\oplus T$. Therefore it is enough to show that $H^0(T)=0$.
Now $T$ is a bundle of rank $k+1$ with $\det T=\OO(-k)$ and thus
$$
T\simeq\wedge^k T^*(-k).
$$
We will prove a little more, indeed we show that
$$
H^0(\wedge^iT^*(-i))=0 \ \ \ \text{for} \ \ \ 1 \ \le i \ \le \ k.
\tag 1.1
$$
We prove (1.1) by induction.
For $i=1$ \ (1.1) holds by the above diagram.
Now, let us consider
$$
0\rig{}E_{|H} \ \rig{} \ \OO_H\oplus T^* \ \rig{} \ \OO_H^k(1) \ \rig{} \ 0.
$$
Let us remember that, since $E$ is a rank $2$ bundle with $c_1(E)=0$, we have
 $\wedge^2E \ \simeq \ \OO$. Therefore, the second wedge power
of $T^*$ twisted by $O_H(-2)$ gives
$$
0 \ \rig{} \ \OO_H(-2) \ \rig{} \ \wedge^2T^*(-2)\oplus T^*(-2) \ \rig{}
\ T^*(-1)^k\oplus\OO_H(-1)^k
$$
which proves (1.1) for $i=2$.
Moreover, the i-th wedge power of $T^*$ twisted by $\OO(-i)$ gives
$$
0 \ \rig{} \ \wedge^iT^*(-i)\oplus \wedge^{i-1}T^*(-i) \ \rig{}
\ \left[\wedge^{i-1}T^*(-i+1)\right]^k\oplus
\left[\wedge^{i-2}T^*(-i+1)\right]^k
$$
and this sequence provides the inductive step.

\qqq
\enddemo

\proclaim {Definition} $W(E)=\{H\in{\PP}^{3*} \ | \ h^0(E_{|H})\neq 0\}$
is called the variety (scheme) of unstable planes of $E$.
Its scheme structure is defined as the degeneracy locus of the mapping
$$
H^1(E(-1))\otimes\OO \ \rig{} \ H^1(E)\otimes\OO(1)
$$
over ${\PP}^{3*}$
(Theorem 1.2 shows that this map drops rank at most by one).
\endproclaim

For an element $A \in \C^{4 \ast} \otimes \C^{k \ast} \otimes \C^{2k+2}$
define the subvariety
$$
X_A \ = \ \{(\overline{f^\ast}, \overline{b^\ast})
\in \PP^{3 \ast} \times \PP^{k-1 \ast} \ \ |
\ \ f^\ast \otimes b^\ast \in \Im(\beta(A, \cdot )) \}.
$$

\proclaim {Lemma 1.3} Let $q_1$ be the projection of 
$\PP^{3 \ast} \times \PP^{k-1 \ast}$ on $\PP^{3 \ast}$.
We have $W(E)=q_1(X_A)$ and the fiber of the projection
$X_A\to q_1(X_A)$ over $H$ is isomorphic to ${\PP}(H^0(E_{|H}))$.
\endproclaim
\demo{Proof} With the notations of the proof of Theorem 1.2 we have that
$H\in W(E)$ iff $h^0(K_{|H}) \neq 0$.
We have $H^0(K_{|H}) = \Ker (\C^{2 k+2}\rig{}\C^4\otimes \C^4)$.
Now suppose that $\overline{f^\ast}$ corresponds to $H$, then the existence 
of a nonzero $\alpha\in H^0(K_{|H})$ is equivalent to
$\beta(A, \alpha)= {f^\ast}\otimes b^\ast$, where
$(\overline{f^\ast} , \overline{b^\ast}) \in \PP^{3 \ast} \times \PP^{k-1 \ast}$ .

\qqq
\enddemo
\proclaim {Corollary 1.4} The morphism $X_A \ \to \ q_1(X_A)$ is an isomorphism
of the underlying varieties, in particular $\dim X_A = \dim q_1(X_A)$.\qqq
\endproclaim

\noindent Recall that special 't Hooft bundles are the instanton bundles
such that $h^0(E(1))=2$. They can be defined through the Serre correspondence
by $k+1$ skew lines lying on a smooth quadric surface [H].
We need the following  special case of a theorem of J. Coanda [Co].
\proclaim {Theorem 1.5 }
If $E$ is an instanton bundle such that $\dim W(E)\ge 2$, then
$E$ is a special 't Hooft bundle and $W(E)$ is a quadric surface. \qqq
\endproclaim

It is known ([H]) that special 't Hooft bundles
are smooth points in the moduli space of instanton bundles.
\proclaim {Corollary 1.6} If $\dim X_A\ge 2$,
then $A$ corresponds to a smooth point in the moduli space of
mathematical instanton bundles.\qqq
\endproclaim

A reformulation of this Corollary into the invariant theoretical
language is as follows.
\proclaim{Corollary 1.7} For any $A^0 \in I$ such that $\dim (X_{A^0}) \ge 2$,
we have $$\dim (T_{\pi (A^0)}MI_{\PP^3}(k)) \ = \ 8k-3.$$
\endproclaim

\proclaim{Lemma 1.8} Suppose
$A^0 \in I$; then $\dim (T_{A^0} I) > 3k^2 + 13k$ if and only if there exists
$0 \neq S^0 \in S^2 \C^4 \otimes \wedge^2 \C^k$ such that
$\xi (A^0, S^0) = 0$, where $\xi$ is the
bilinear $\SL_4 \times \SL_k \times \SP_{2k+2}$-morphism defined by
$$
\gathered
\xi : \C^{4 \ast} \otimes \C^{k \ast} \otimes \C^{2k+2}
\ \times \ S^2 \C^4 \otimes \wedge^2 \C^k
\ \rightarrow \ \C^4 \otimes \C^k \otimes \C^{2k+2}, \\
(f_i^\ast \otimes b_j^\ast \otimes h_l,
\ f_{i'} f_{i''} \otimes b_{j'} \wedge b_{j''})
\ \mapsto \ (\delta_{ii'} f_{i''} + \delta_{i''i} f_{i'})
\otimes (\delta_{j'j} b_{j''} - \delta_{j''j} b_{j'}) \otimes h_l.
\endgathered
$$

\endproclaim
\demo{Proof} We have $\dim (T_{A^0} I) > 3k^2 + 13k$ iff the
differential \ $d \gamma|_{A^0}$ \ is nonsurjective.
The differential $d \gamma|_{A^0}$ is nonsurjective iff
$(d \gamma|_{A^0})^\ast$ is noninjective, i.e.,
$(d \gamma|_{A^0})^\ast(S^0) = 0$ for some
$0 \neq S^0 \in S^2 \C^4 \otimes \wedge^2 \C^k$.
It can be easily checked that
$$
(d \gamma|_A)^\ast(S) \equiv \xi (A,S).
$$
Hence, $\dim (T_{A^0} I) > 3k^2 + 13k$ iff $\xi (A^0,S^0) = 0$ for some
$0 \neq S^0 \in S^2 \C^4 \otimes \wedge^2 \C^k$.

\qqq
\enddemo

 For the convenience of the reader we give a cohomological
interpretation of Lemma 1.8. Let $E^0$ be the instanton bundle defined by
$A^0\in I$ as the cohomology bundle of monad (1.1). By Lemma 1.1 and
deformation theory the assumption $\dim (T_{A^0} I) > 3k^2 + 13k$  is equivalent to
$h^1(S^2E^0)=\dim (T_{\pi (A^0)}MI_{\PP^3}(k))>8k-3$. Therefore, the assumption
of Lemma 1.8 is equivalent to $H^2(S^2E^0)\neq 0$. The second symmetric power
of the left hand side of (1.1) gives $H^2(S^2E^0)\simeq H^2(S^2(\Ker~f_{A^0}))$.
The second symmetric power of the right hand side of (1.1) gives
$$
H^2 (S^2(\Ker~f_{A^0}))\simeq
\Coker \left[ H^0(\OO(1))\otimes \C^{k \ast} \otimes \C^{2k+2 \ast} \ \rig{\Phi}
\ H^0(\OO(2)) \otimes \wedge^2 (\C^{k \ast})\right].
$$
Lemma 1.8 follows because the dual of $\Phi$ can be identified with
$\xi (A^0,\cdot)$.


\head \S2. Algebraic lemmas. \endhead

Among all this section we prove some algebraic lemmas that we will use
in order to prove our main result.

\proclaim{Lemma 2.1}
Suppose $R$ is a nonzero block-matrix:
$$
R = \left( \matrix R^1 \\ R^2 \endmatrix \right),
$$
where $R^i$ is a skew-symmetric matrix of size $k \times k$,
\ then there exists a column $v_0$ of height $k$ such that
$$
R v_0 \ = \ \left( \matrix
\lambda_1 u_0 \\ \lambda_2 u_0
\endmatrix \right) \ \neq 0
$$
for some column $u_0$ of height $k$, \ $\lambda_1, \lambda_2 \in \C$.
\endproclaim
\demo{Proof} Suppose that $det(R^1) \neq 0$. In this case set
$v_0 \in \Ker (R^2 - \mu_0 R^1)$, where $\mu_0$ is a root
of the equation $det(R^2 - \mu R^1) = 0$.

Suppose that $det(R^1) = 0$. One can assume that
$$
R^1 \ = \ \left( \matrix R^1_{11} & 0 \\ 0 & 0
\endmatrix \right), \ \ \ \ \
R^2 = \left( \matrix R^2_{11} & R^2_{12} \\ R^2_{21} & R^2_{22}
\endmatrix \right),
$$
where $R^1_{11}$ is a skew-symmetric matrix of size
$k' \times k'$, \ $k' < k$, \ $det (R^1_{11}) \neq 0$ and $R^2_{11}$
is a skew-symmetric matrix of size $k' \times k'$.
If $R^2_{12} \neq 0$ or $R^2_{22} \neq 0$, then we set
$v_0 = \left( \matrix 0 \\ v_0' \endmatrix \right)$ for some $v_0'$ such that
$R^2_{12} v'_0 \neq 0$ or $R^2_{22} v'_0 \neq 0$.
If $R^2_{12} = 0$ and $R^2_{22} = 0$, then $R^2_{21} = 0$ and we set
$v_0 = \left( \matrix v_0' \\ 0 \endmatrix \right)$, where
$$
\left( \matrix R^1_{11} \\ R^2_{11} \endmatrix \right) v_0'
\ = \ \left( \matrix
\lambda_1 u_0' \\ \lambda_2 u_0'
\endmatrix \right) \ \neq 0.
$$
\qqq
\enddemo

Consider the linear spaces $\C^4$ and $\C^k$. Let $f_1, ... , f_4$
be the standard basis of $\C^4$ and let $f_1^\ast, ... , f_4^\ast$ be
the dual basis of the dual space $\C^{4 \ast}$. Let $b_1, ... , b_k$
be the standard basis of $\C^k$ and let $b_1^\ast, ... , b_k^\ast$ be
the dual basis of the dual space $\C^{k \ast}$. The group $\SL_4$
acts canonically on the space $\C^4$ and the group $\SL_k$
acts canonically on the space $\C^k$. So the actions of the
group $\SL_4 \times \SL_k$ are defined on the spaces
$\C^4, \ \C^{4 \ast}, \ \C^k, \ \C^{k \ast},
\ \C^4 \otimes \C^k, ... $.

Consider the linear space $S^2 \C^4 \otimes \wedge^2 \C^k$.
For an element $S \in S^2 \C^4 \otimes \wedge^2 \C^k$ define
$$
\rk (S) = \dim (\Im (\rho (S, \cdot))),
$$
where
$$
\gathered
\rho : S^2 \C^4 \otimes \wedge^2 \C^k
\ \times \ \C^{4 \ast} \otimes \C^{k \ast}
\ \rightarrow \ \C^4 \otimes \C^k, \\
(f_{i'} f_{i''} \ \otimes \ b_{j'} \wedge b_{j''},
\ f_i^\ast \otimes b_j^\ast) \ \mapsto
\ (\delta_{i'i} f_{i''} + \delta_{i''i} f_{i'}) \otimes
(\delta_{j'j} b_{j''} - \delta_{j''j} b_{j'})
\endgathered
$$
is the bilinear $\left( \SL_4 \times \SL_k\right)$-morphism. Note that
$\rk (S)$ is an even number. The following lemma is the only place in the paper 
where we need the assumption $k\le 5$.

\proclaim{Lemma 2.2}
Suppose $2 \le k \le 5$ and consider \ $S \in S^2 \C^4 \otimes \wedge^2 \C^k$,
such that $2 \le \rk (S) \le 2k-2$. Then one of the following
conditions holds:
\roster
\item
$\rho (S, B^{\ast 0}) = f^0 \otimes b^0 \neq 0$ for some
$B^{\ast 0} \in \C^{4 \ast} \otimes \C^{k \ast}$, \ $f^0 \in \C^4$,
\ $b^0 \in \C^k$.
\item $\rk (S) = 6$ and there exists \ $0 \neq f^{\ast 0} \in \C^{4 \ast}$
\ such that $\rho (S, f^{\ast 0} \otimes b^\ast) = 0$
for all $b^\ast \in \C^{k \ast}$.
\item $\rk (S) = 8$ and \ $\dim (Z_S) \ge 2$, where
$$
\gathered
Z_S \ = \ \{(\overline {f^\ast}, \overline {b^\ast})
\in \PP^{3 \ast} \times \PP^{k-1 \ast} \ \ |
\ \ \rho (S, f^\ast \otimes b^\ast) = 0 \}, \\
\ \ \PP^{3 \ast} \ = \ P \C^{4 \ast},
\ \ \PP^{k-1 \ast} \ = \ P \C^{k \ast}.
\endgathered
$$
\endroster
\endproclaim
\demo{Proof}
Consider the coordinate expression of $S$
in the bases $\{ f_i \}$ and $\{ b_i \}$:
$$
S = \sigma_{lp}^{ij} f_l f_p \otimes b_i \wedge b_j.
$$
We get a block matrix $\sigma$ defined by
$$
\sigma \ = \ (\sigma^{ij})_{1 \le i,j \le k} \ =
\ \left( \matrix
     0      &   \sigma^{12} & \hdots & \sigma^{1k} \\
\sigma^{21} &         0     & \hdots & \sigma^{2k} \\
\vdots      &     \vdots    & \ddots & \vdots \\
\sigma^{k1} &  \sigma^{k2}  & \hdots &    0
\endmatrix \right),
$$
where $\sigma^{ij} = (\sigma_{lp}^{ij})_{1 \le l,p \le 4}$
is a symmetric matrix of size $4 \times 4$,
\ $\sigma^{ij} = - \sigma^{ji}$.
There is a second coordinate expression
$$
S = \hsigma_{lp}^{ij} f_i f_j \otimes b_l \wedge b_p,
$$

and we get a second block matrix $\hsigma$ defined by
$$
\hsigma \ = \ (\hsigma^{ij})_{1 \le i,j \le 4} \ =
\ \left( \matrix
\hsigma^{11} & \hsigma^{12} & \hsigma^{13} & \hsigma^{14} \\
\hsigma^{21} & \hsigma^{22} & \hsigma^{23} & \hsigma^{24} \\
\hsigma^{31} & \hsigma^{32} & \hsigma^{33} & \hsigma^{34} \\
\hsigma^{41} & \hsigma^{42} & \hsigma^{43} & \hsigma^{44}
\endmatrix \right),
$$
where $\hsigma^{ij} = (\hsigma_{lp}^{ij})_{1 \le l,p \le k}$
is a skew-symmetric matrix of size $k \times k$,
\ $\hsigma^{ij} = \hsigma^{ji}$.

Transform the basis $\{ b_i \}$ and obtain
$$
r \ \overset{\text{def}}\to{=} \ \rk(\sigma^{12}) \ =
\ \underset {\{(c_{ij})_{1 \le i,j \le k}\}}\to {max}
\{ \rk(c_{1i} \sigma^{ij} c_{2j}) \}. \tag 2.1
$$
We have
$$
2k-2 \ \ge \ \rk (S) \ = \ \rk (\sigma) \ =
\ \rk (\hsigma) \ge \ 2\rk (\sigma^{12}) \ = \ 2r.
$$
Therefore, one of the following cases holds:
\roster
\item"(a)" $r = 1$ or $2$,
\item"(b)" $r = 3$, $\rk (\sigma) = 6$, and $k \ge 4$,
\item"(c)" $r = 4$, $\rk (\sigma) = 8$, and $k = 5$,
\item"(d)" $r = 3$, $\rk (\sigma) = 8$, and $k = 5$.
\endroster

Transform the basis $\{ f_i \}$ and obtain
$$
\sigma_{lp}^{12} \ = \cases
1 \ \ \ \ \ \text{if} \ \ \ 1 \le l = p \le r, \\
0 \ \ \ \ \ \text{if} \ \ \ l \neq p \ \ or \ \ l = p > r.
\endcases \tag 2.2
$$
From (2.1) it follows that
$\sigma_{lp}^{ij} = 0$ for $l,p > r$ whence
$$
\hsigma^{ij} = 0 \ \ \ \ \ \text{for} \ \ i,j > r. \tag 2.3
$$

\bf (a). \rm Consider the case (a).

In this case we prove that the condition (1) holds, i.e.,
we prove that there exists a column $f^0$ of height $4$
and columns
$b^0, B^{\ast 0 1} , \dots , B^{\ast 0 4}$ of height $k$ such that
$$
\hsigma \left( \matrix
B^{\ast 0 1} \\ \vdots \\ B^{\ast 0 4}
\endmatrix \right)
\ = \ \left( \matrix
f^0_1 b^0 \\ \vdots \\ f^0_4 b^0
\endmatrix \right)
\ \neq \ 0.
$$
But this easily follows from Lemma 2.1 and (2.3).

\bf (b). \rm Consider the case (b).

In this case we prove that the condition (2) holds, i.e.
we prove that there exists a column $f^{\ast 0}$ of
height $4$ such that
$$
\sigma \left( \matrix
b^\ast_1 f^{\ast 0} \\ \vdots \\ b^\ast_k f^{\ast 0}
\endmatrix \right) \ = \ 0 \tag 2.4
$$
for any column $b^\ast$ of height $k$.
From the condition $\rk (\sigma) = 6$ and (2.1) it follows that
$$
\sigma^{ij} \ = \ \left( \matrix
\sigma^{ij}_{11} & \sigma^{ij}_{12} & \sigma^{ij}_{13} & 0 \\
\sigma^{ij}_{21} & \sigma^{ij}_{22} & \sigma^{ij}_{23} & 0 \\
\sigma^{ij}_{31} & \sigma^{ij}_{32} & \sigma^{ij}_{33} & 0 \\
       0         &        0         &        0         & 0
\endmatrix \right).
$$
From this for
$$
f^{\ast 0} \ = \ \left( \matrix
0 \\ 0 \\ 0 \\ 1
\endmatrix \right)
$$
it easily follows (2.4)

\bf (c). \rm Consider the case (c).

In this case we prove that the condition (3) holds. We have:
$$
Z_S \ =
\ \{ (\overline{f^\ast}, \ \overline{b^\ast})
\ = \ ( \overline{\left( \matrix
f^\ast_1 \\ \vdots \\ f^\ast_4
\endmatrix \right)},
\ \overline{\left( \matrix
b^\ast_1 \\ \vdots \\ b^\ast_5
\endmatrix \right)} )
\ \ | \ \ \sigma \left( \matrix
b^\ast_1 f^\ast \\ \vdots \\ b^\ast_5 f^\ast
\endmatrix \right) = 0 \}.
$$

Consider the matrix
$$
\tsigma \ =
\ \left( \matrix
     0      &    \EE_4    & \sigma^{13} & \sigma^{14} & \sigma^{15} \\
  -\EE_4    &      0      & \sigma^{23} & \sigma^{24} & \sigma^{25}
\endmatrix \right),
$$
where $\EE_4$ is the identity matrix of size $4 \times 4$. The $8$ rows
of the matrix $\tsigma$ are the first $8$ rows of the matrix $\sigma$.
Since $\rk (\sigma) = 8 = \rk (\tsigma)$, for a matrix $P$ of
size $20 \times p$ we have:
$$
\sigma P = 0 \ \ \ \ \ \ \ \text{iff}
\ \ \ \ \ \ \ \tsigma P = 0. \tag 2.5
$$

For $3 \le i \le 5$ consider the following matrix $P_i$ of size
$20 \times 4$:
$$
P_i = \ \left( \matrix
-\sigma^{2i} \\ \sigma^{1i} \\ P_{i3} \\ P_{i4} \\ P_{i5}
\endmatrix \right),
$$
where $P_{ii} = -\EE_4$ and $P_{ij} = 0$ for $j \neq i$. We see that
$\tsigma \cdot P_i = 0$.
From (2.5) it follows that $\sigma \cdot P_i = 0$ or
$$
\sigma^{ji} = \sigma^{1j} \sigma^{2i} - \sigma^{2j} \sigma^{1i},
\ \ \ \ \ 3 \le j \le 5.
$$
From this we obtain
$$
\split
0 \ = & \  \sigma^{ji} + (\sigma^{ij})^\top =
\sigma^{1j} \sigma^{2i} \ - \ \sigma^{2j} \sigma^{1i} \ +
\ (\sigma^{1i} \sigma^{2j} \ - \ \sigma^{2i} \sigma^{1j})^\top \ = \\
& [\sigma^{1j}, \sigma^{2i}] \ + \ [\sigma^{1i}, \sigma^{2j}],
\ \ \ \ \ \ \ \ 3 \le i,j \le 5.
\endsplit
$$
One can rewrite these equations into the following compact form:
$$
[ t_1 \sigma^{13} + t_2 \sigma^{14} + t_3 \sigma^{15},
\ t_1 \sigma^{23} + t_2 \sigma^{24} + t_3 \sigma^{25} ] \ = \ 0 \tag 2.6
$$
for all $t_1, t_2, t_3 \in \C$.

\proclaim{Claim 1} For every $(b_3^\ast, b_4^\ast, b_5^\ast) \neq (0,0,0)$
there exists $(b_1^\ast, b_2^\ast)$ and a nonzero column $f^\ast$
of height $4$ such that
$$
\sigma \left( \matrix
b^\ast_1 f^\ast \\ \vdots \\ b^\ast_5 f^\ast
\endmatrix \right) = 0.
$$
\endproclaim
\demo{Proof of Claim 1} From (2.6) it follows
that the symmetric matrices
$$
b^\ast_3 \sigma^{13} + b^\ast_4 \sigma^{14} + b^\ast_5 \sigma^{15},
\ \ \ \ \ b^\ast_3 \sigma^{23} + b^\ast_4 \sigma^{24} + b^\ast_5 \sigma^{25}
$$
commute therefore they have a common eigenvector $f^\ast$ with the
eigenvalues \ $b^\ast_2$, $-b^\ast_1$ \ respectively. We have
$$
\tsigma \left( \matrix
b^\ast_1 f^\ast \\ \vdots \\ b^\ast_5 f^\ast
\endmatrix \right) = 0
$$
and from this and (2.5) Claim 1 \ follows.

\qqq
\enddemo
From Claim 1 \ it follows that $\dim (Z_S) \ge 2$.

\bf (d). \rm Consider the case (d).

In this case we prove that the condition (3) holds, i. e.,
we prove that $\dim (Z_S) \ge 2$.

\proclaim{Claim 2} Suppose $N \subset P\C^{5 \ast}$ is a line
in general position; then there exists $0 \neq f^{\ast 0} \in \C^{4 \ast}$,
$\overline{b^{\ast 0}} \in N$ such that
$\rho (S, f^{\ast 0} \otimes b^{\ast 0}) = 0$.
\endproclaim
\demo{Proof of Claim 2} One can assume that
$N = \overline{\langle b_1, b_2 \rangle}$. We have
to prove that there exists a column $f^{\ast 0}$ of height $4$ and
$b^\ast_1, b^\ast_2 \in \C$,
$(b^\ast_1, b^\ast_2) \neq (0,0)$ such that
$$
\sigma \left( \matrix
b^\ast_1 f^{\ast 0} \\ b^\ast_2 f^{\ast 0} \\ 0 \\ 0 \\ 0
\endmatrix \right) \ = \ 0. \tag 2.7
$$

Consider the $4^{th}$ and $8^{th}$ rows of the matrix $\sigma$:
$$
\gathered
\row_4(\sigma) \ = \ (0, \dots ,0,
\sigma^{13}_{41},\sigma^{13}_{42},\dots,
\sigma^{15}_{43},\sigma^{15}_{44}), \\
\row_8(\sigma) \ = \ (0, \dots,0,
\sigma^{23}_{41},\sigma^{23}_{42},\dots,
\sigma^{25}_{43},\sigma^{25}_{44}).
\endgathered
$$

We want to show that  $\row_4(\sigma)$ and $\row_8(\sigma)$  are linearly 
dependent.
Suppose that $\row_4(\sigma)$ and $\row_8(\sigma)$ are linearly
independent, then the first $8$ rows of the matrix $\sigma$
are linearly independent. Since $\rk (\sigma) = 8$, we see that every
row of $\sigma$ is a linear combination of the first $8$ rows.
From $\row_4(\sigma) \neq 0$ it follows that
$\sigma^{1i}_{4j} \neq 0$ for some $3 \le i \le 5$, $1 \le j \le 4$.
Since $\sigma^{i1}_{j4} = -\sigma^{1i}_{4j} \neq 0$, we see that
$(4(i-1)+j)$th row
$$
\row_{4(i-1)+j}(\sigma) =
(\sigma^{i1}_{j1},\sigma^{i1}_{j2},
\sigma^{i1}_{j3},\sigma^{i1}_{j4},
\sigma^{i2}_{j1},\sigma^{i2}_{j2},
\sigma^{i2}_{j3},\sigma^{i2}_{j4}, \dots )
$$
of the matrix $\sigma$ is \it not \rm a linear combination
of the first $8$ rows. This contradiction proves that
$\row_4(\sigma)$ and $\row_8(\sigma)$ are linearly dependent.

Finally, to obtain (2.7) we take
$$
f^{\ast 0} \ = \ \left( \matrix
0 \\ 0 \\ 0 \\ 1
\endmatrix \right),
$$
and $b^\ast_1, b^\ast_2$ \ such that \ $(b_1^\ast, b_2^\ast) \neq (0,0)$ \ and
\ $b^\ast_1 \row_4(\sigma) + b_2 \row^\ast_8(\sigma) = 0$.

\qqq
\enddemo
From Claim 2 \ it follows that $\dim (Z_S) \ge 3 > 2$.

\qqq
\enddemo


\head \S3. The proof of Theorem 0.1.
\endhead

\proclaim{Lemma 3.1} Consider elements
$A^0 \in \C^{4 \ast} \otimes \C^{k \ast} \otimes \C^{2k+2}$ and
\ $S^0 \in S^2 \C^4 \otimes \wedge^2 \C^k$ such that
$\xi (A^0, S^0) = 0$.
\roster
\item Suppose
$$
\tau : \C^{4 \ast} \otimes \C^{k \ast} \otimes \C^{2k+2}
\ \times \ S^2 \C^4 \otimes \wedge^2 \C^k
\ \times \ \C^{2k+2}
\ \rightarrow \ \C^4 \otimes \C^k
$$
is an arbitrary trilinear
$\left(\SL_4 \times \SL_k \times \SP_{2k+2}\right)$-morphism;
then $$\tau (A^0, S^0, \C^{2k+2}) = 0$$
\item Suppose
$$
\alpha : \C^{4 \ast} \otimes \C^{k \ast} \otimes \C^{2k+2}
\ \times \ S^2 \C^4 \otimes \wedge^2 \C^k
\ \times \ \C^{4 \ast} \otimes \C^{k \ast}
\ \rightarrow \ \C^{2k+2}
$$
is an arbitrary trilinear
$\left(\SL_4 \times \SL_k \times \SP_{2k+2}\right)$-morphism;
then $$\alpha (A^0, S^0, \C^{4 \ast} \otimes \C^{k \ast}) = 0$$
\endroster
\endproclaim
\demo{Proof} 

(1) Consider the following nontrivial trilinear
$\SL_4 \times \SL_k \times \SP_{2k+2}$-morphism:
$$
\gathered
\tau_0 : \C^{4 \ast} \otimes \C^{k \ast} \otimes \C^{2k+2}
\ \times \ S^2 \C^4 \otimes \wedge^2 \C^k
\ \times \ \C^{2k+2}
\ \rightarrow \ \C^4 \otimes \C^k, \\
(A, S, h) \ \mapsto \ \kappa(\xi (A, S), h),
\endgathered
$$
where
$$
\gathered
\kappa : \C^{4 \ast} \otimes \C^{k \ast} \otimes \C^{2k+2}
\ \times \ \C^{2k+2}
\ \rightarrow \ \C^4 \otimes \C^k, \\
(f_i \otimes b_j \otimes h_l, h_{l'})
\ \mapsto \ \omega (h_l, h_{l'}) f_i \otimes b_j
\endgathered
$$
is the bilinear
$\SL_4 \times \SL_k \times \SP_{2k+2}$-morphism.

On the other hand the $\left(\SL_4 \times \SL_k \times \SP_{2k+2}\right)$-module
$$
(\C^{4 \ast} \otimes \C^{k \ast} \otimes \C^{2k+2})
\ \otimes \ (S^2 \C^4 \otimes \wedge^2 \C^k)
\ \otimes \ \C^{2k+2}
$$
contains the irreducible
$\left(\SL_4 \times \SL_k \times \SP_{2k+2}\right)$-module
\ $\C^4 \otimes \C^k$ \ with multiplicity 1. \ Therefore, there
exists a unique, up to a scalar factor, nontrivial trilinear
$\left(\SL_4 \times \SL_k \times \SP_{2k+2}\right)$-morphism
$$
\C^{4 \ast} \otimes \C^{k \ast} \otimes \C^{2k+2}
\ \times \ (S^2 \C^4 \otimes \wedge^2 \C^k)
\ \times \ \C^{2k+2}
\ \rightarrow \ \C^4 \otimes \C^k.
$$
Thus, $\tau = c \tau_0$ for some $c \in \C$ and we get
$$
\tau (A^0, S^0, \C^{2k+2}) \ = \ c \tau_0 (A^0, S^0, \C^{2k+2})
\ = \ c \kappa (\xi (A^0, S^0), \C^{2k+2}) \ = \ 0.
$$

(2) Consider the linear mapping
$\alpha (A^0, S^0, \cdot)^\ast$
dual to $\alpha (A^0, S^0, \cdot)$ . From (1) it follows that
$\alpha (A^0, S^0, \cdot)^\ast \ = \ 0$. Thus,
$\alpha (A^0, S^0, \C^{4 \ast} \otimes \C^{k \ast}) \ = \ 0$.

\qqq
\enddemo

\demo{Proof of Theorem 0.1} 

We suppose that there exists $A^0 \in I$ such that
$\dim (T_{\pi (A^0)}MI_{\PP^3}(k)) > 8k-3$, \ \ $2 \le k \le 5$ and we 
obtain a contradiction.

From Corollary 1.7 it follows that
$$
\dim (X_{A^0}) \ \le \ 2 \tag 3.1
$$
and by Lemma 1.1 we have $\dim (T_{A^0} I) > 3k^2 + 18k$.
Hence, by Lemma 1.8  there exists
$0 \neq S^0 \in S^2 \C^4 \otimes \wedge^2 \C^k$ such that
$$
\xi (A^0, S^0) = 0. \tag 3.2
$$
Consider the following composition of linear mappings
$$
\rho (S^0, \cdot) \circ \beta(A^0, \cdot) \ :
\ \C^{2k+2} \rightarrow \C^4 \otimes \C^k,
\ \ \ \ \ h \ \mapsto \ \rho (S^0, \beta (A^0, h)),
$$
where $\beta$ is defined in \S 1 and $\rho$ is defined in \S 2.
From (3.2) and Lemma 1.8 (1) it follows that
$\rho (S^0, \cdot) \circ \beta(A^0, \cdot) = 0$ or
$$
\Im(\beta(A^0, \cdot)) \subset \Ker(\rho(S^0, \cdot)). \tag 3.3
$$
On the other hand, by (E3) we have $\rk(\beta (A^0, \cdot)) = 2k+2$ and
with (3.3) this gives us
$$
\rk(\rho (S^0, \cdot)) \le 2k-2. \tag 3.4
$$

Therefore, from (3.4) it follows that one of the conditions (1) - (3) of Lemma 2.2
holds for $S = S^0$.

\bf I. \rm Consider the case when the condition (1) of Lemma 2.2
holds for $S = S^0$.

In this case, consider the following composition of 
linear mappings
$$
\varepsilon (A^0, \cdot) \circ \rho (S^0, \cdot) \ :
\ \C^{4 \ast} \otimes \C^{k \ast}
\ \rightarrow \ \C^{2k+2},
\ \ \ \ \ B^\ast \ \mapsto \ \varepsilon (A^0, \rho (S^0, B^\ast)).
$$
From Lemma 1.8 (2) it follows that
$\varepsilon (A^0, \cdot) \circ \rho (S^0, \cdot)) = 0$.
By the condition (1) of Lemma 2.2 there exists
$B^{\ast 0} \ \in \ \C^{4 \ast} \otimes \C^{k \ast}$ such
that $\rho (S^0, B^{\ast 0}) = f^0 \otimes b^0 \neq 0$. Thus, we have
$\varepsilon (A^0, f^0 \otimes b^0) = \varepsilon (A^0, \rho (S^0, B^{\ast 0})) = 0$
and therefore $A^0 \notin I_1$. But this contradicts the fact that $A^0 \in I$.

\bf II. \rm Consider the case when the condition (2) of Lemma 2.2
holds for $S = S^0$.

From (3.4) it follows that $k = 4$ or $k = 5$.
By the condition (2) of Lemma 2.2 we have
$\{ f^{\ast 0} \} \times \C^{k \ast} \subset \Ker(\rho(S^0, \cdot))$.
On the other hand, we have (3.3) and
$$
\dim (\Ker(\rho(S^0, \cdot))) - \dim (\Im(\beta(A^0, \cdot))) \ =
\ \cases 0 \ \ \ \ \ \text{if} \ k=4, \\
         2 \ \ \ \ \ \text{if} \ k=5.
\endcases
$$
Therefore
$\Im(\beta(A^0, \cdot) \supset \{ f^{\ast 0} \} \times M$ for some linear
subspace $M \subset \C^{k \ast}$ of dimension $\ge 3$.
But this contradicts (3.1).

\bf III. \rm Consider the case when the condition (3) of Lemma 2.2
holds for $S = S^0$.

From (3.4) it follows that $k = 5$.
Thus $$\dim (\Im(\beta(A^0, \cdot))) = 12 = \dim (\Ker(\rho(S^0, \cdot)))$$
and from this, together with (3.3), it follows that
$\Im(\beta(A^0, \cdot)) = \Ker(\rho(S^0, \cdot))$. Therefore
$X_{A^0} = Z_{S^0}$. From this and the condition (3) of Lemma 2.2 we
obtain $\dim (X_{A^0}) = \dim (Z_{S^0}) \ge 2$. But this again
contradicts (3.1).

\qqq
\enddemo

\enddocument